\documentclass[12pt,a4paper,oneside]{amsart}
\usepackage{amsfonts, amsmath, amssymb, amsthm, amscd, hyperref}
\usepackage[T2A]{fontenc}
\usepackage[cp1251]{inputenc}
\usepackage[english]{babel}
\usepackage[dvips]{graphicx}
\usepackage{subfig}
\usepackage{psfrag}
\usepackage{epstopdf}
\usepackage{cmap}
\usepackage{enumerate}
\usepackage[english]{babel}

\newtheorem{lemma}{Lemma}
\newtheorem{remark}{Remark}
\newtheorem{theorem}{Theorem}
\newtheorem{corollary}{Corollary}
\newtheorem{conjecture}{Conjecture}
\newtheorem{question}{Question}

\newcommand{\la}{\lambda}

\newcommand{\e}{\varepsilon}

\def\frl{\forall}%
\newcommand{\men}{\leqslant}
\newcommand{\bol}{\geqslant}
\newcommand{\bra}{\langle}
\newcommand{\ket}{\rangle}
\newcommand{\R}{\mathbb{R}}

\def\B{\mathfrak{B}}%
\def\SS{\partial\B_1(o)}%
\def\SSS{\partial\B_1^*(o)}%
\def\BB{\B_1(o)}%
\newcommand{\reff}[1]{(\ref{#1})}%

\newcommand{\norm}[1]{\left\| #1 \right\|}

\def\vn{\mathop{\rm int}}

\newcommand{\mglx}[1]{\rho_{X}\!\!\!\:\left( #1 \right)}\!%
\newcommand{\mcox}[1]{\delta_{X}\!\left( #1 \right)}%
\newcommand{\mcoxt}[2]{\delta_{X}\!\left( #1, #2 \right)}%
\newcommand{\mgbx}[1]{\delta_{X}^{+}\!\left( #1 \right)}%
\newcommand{\mgbxt}[2]{\delta_{X}^{+}\!\left( #1, #2 \right)}%
\newcommand{\lamx}[1]{\lambda^{-}_{X}\!\left( #1 \right)}%
\newcommand{\lapx}[1]{\lambda^{+}_{X}\!\left( #1 \right)}%
\newcommand{\phipx}[1]{\varphi^{+}_{X}\!\left( #1 \right)}%
\newcommand{\phimx}[1]{\varphi^{-}_{X}\!\left( #1 \right)}%

\newcommand{\zetap}[1]{\zeta^{+}_{X}\!\left( #1 \right)}%
\newcommand{\zetam}[1]{\zeta^{-}_{X}\!\left( #1 \right)}%

\newcommand{\gammap}[1]{\gamma^{+}_{X}\!\left( #1 \right)}%
\newcommand{\gammam}[1]{\gamma^{-}_{X}\!\left( #1 \right)}%

\newcommand{\prp}{\urcorner}%
\newcommand{\prooff}{\noindent {\bf Proof.}\\}%
\newcommand{\bbox}{\par\noindent\ensuremath{\Box}\par\noindent}%
\newenvironment{prf}
{\prooff}
{\bbox}

\begin{document}

\begin{center}
\textbf{\Large New Moduli for Banach Spaces}\\[0.3cm]

\emph{Grigory Ivanov\footnote{
DCG, FSB, Ecole Polytechnique F\'ed\'erale de Lausanne, Route Cantonale, 1015 Lausanne, Switzerland.
\\ \noindent
Department of Higher Mathematics, Moscow Institute of Physics and Technology,  Institutskii pereulok 9, Dolgoprudny, Moscow
region, 141700, Russia.
\\ \noindent
grimivanov@gmail.com
\\ \noindent 
Research partially supported by Swiss National Science Foundation grants 200020-165977 and 200021-162884
Supported by Russian Foundation for Basic Research, project 16-01-00259.
} 
and Horst Martini\footnote{Faculty of Mathematics , TU Chemnitz, 09107 Chemnitz, Germany}}

\end{center}
\vspace{0.2cm}

{\bf Abstract.}
Modifying the moduli of supporting convexity and supporting smoothness, we introduce new moduli for Banach spaces which occur, e.g., as lengths of catheti of right-angled triangles (defined via so-called quasi-orthogonality). These triangles have two boundary points of the unit ball  of a Banach space as endpoints of their hypotenuse, and their third vertex lies in a supporting hyperplane of one of the two other vertices. Among other things it is our goal to quantify via such triangles the local deviation of the unit sphere from its supporting hyperplanes. We prove respective Day-Nordlander type results, involving generalizations of the modulus of convexity and the modulus of Bana\'{s}. 

\bigskip


Mathematics Subject Classification (2010): 46B07, 46B20, 52A10, 52A20, 52A21
\bigskip

Keywords: Birkhoff-James orthogonality, Day-Nordlander type results, Milman modulus, modulus of (supporting) convexity,
modulus of (supporting) smoothness, quasi-orthogonality

\section{Introduction}
The modulus of convexity (going back to \cite{Clarkson_ucsp}) 
and the modulus of smoothness (defined in \cite{Day_uncsp})
are well known classical constants from Banach space theory. 
For these two notions various interesting applications were found, 
and a large variety of natural refinements, generalizations, and modifications of
them created an impressive bunch of interesting results and problems; 
see, e.g., \cite{Day-Nord}, \cite{Milman_eng}, \cite{Benitez_Rectangular_constant}, \cite{Banass1},
 \cite{onmoduli_Banas}, and \cite{Ivanov_supp_modulus}, to cite only references close to our discussion here. Inspired by \cite{Banass1}, two further constants in this direction were introduced and investigated in \cite{Ivanov_supp_modulus}, namely the modulus of supporting convexity and the modulus of supporting smoothness. These moduli suitably quantify the local deviation of the boundary of the unit ball of a real Banach space from its supporting hyperplanes near to arbitrarily chosen touching points. Using the concept of right-angled triangles in terms of so-called quasi-orthogonality (which is closely related to the concept of Birkhoff-James orthogonality), 
we modify and complete the framework of moduli defined in \cite{Clarkson_ucsp}, \cite{Banass1}, 
and \cite{Ivanov_supp_modulus} by introducing and studying new related constants. These occur as lengths of catheti of such triangles, whose hypotenuse connects two boundary points of the unit ball and whose third vertex lies in the related supporting hyperplane. We prove Day-Nordlander type results referring to these moduli, yielding even generalizations of the constants introduced in \cite{Clarkson_ucsp}, \cite{Banass1}, 
and \cite{Ivanov_supp_modulus}. Respective results on Hilbert spaces are obtained, too. At the end we discuss some conjectures and questions which refer to further related inequalities between such moduli (for general Banach spaces, but also for Hilbert spaces),  possible characterizations of inner product spaces, and  Milman's moduli.

The paper is organized as follows: After presenting our notation and basic definitions in Section \ref{section_def}, we clarify the geometric position of the mentioned right-angled triangles close to a point of the unit sphere of a Banach space and its corresponding supporting hyperplane. This yields a clear geometric presentation of the new moduli, but also of further moduli already discussed in the literature. In Section \ref{section_catheti} we particularly study properties of the catheti of these triangles, yielding also the announced results of Day-Nordlander type and results on Hilbert spaces. In a similar way, we study properties of the hypotenuses in Section \ref{section_hypothenuse}, obtaining again Day-Nordlander type results and further new geometric inequalities. In Section \ref{section_monotonicity} our notions and results are put into a more general framework, connected with concepts like monotone operators, dual mappings of unit spheres and their monotonicity. And in Section \ref{section_open_questions} some open questions and conjectures on the topics shortly described above are collected.
\section{Notation and basic  definitions}\label{section_def}
In the sequel we shall need the following notation.
Let $X$ be a {\it real Banach space}, and $X^*$ be its {\it conjugate space}. 
We use  $H$ to denote a {\it Hilbert space}. For a set $A \subset X$ we denote by $  \partial A$ and $\vn A$ the \emph{boundary} and the \emph{interior} of  $A,$ respectively.
We use $\bra p,x \ket$ to denote the \emph{value of a functional} $p \in X^*$ \emph{at a vector} $x \in X.$
For $R>0$ and $c \in X$ we denote by $\B_R(c)$ the \emph{closed ball} 
{\it with center $c$ and radius $R,$} and
by $\B_R^{*}(c)$ the respective \emph{ball in the conjugate space}. Thus, $\partial \B_1(o)$ denotes the \emph{unit sphere} of $X$.
By definition, we put $J_1(x) = \{p \in \SSS :\, \bra p, x \ket\ = \norm{x}\}.$

We will use the notation $xy$ for the {\it segment} with the (distinct) endpoints $x$ and $y$, for the {\it line} passing through these points, for (oriented) {\it arcs} from $\partial \B_R(c)$, as well as for the {\it vector} from $x$ to $y$ (the respective meaning will always be clear by the context). Further on, abbreviations like $abc$ and $abcd$ are used for {\it triangles} and {\it 4-gons} as convex hulls of these three or four points.

We say that $y$ is {\it quasi-orthogonal}  to the vector $x \in X\setminus \{o\}$
and write  $y\urcorner x$
if there exists a functional $p \in J_1(x)$  such that  $\bra p, y \ket = 0.$
Note that  the following conditions are equivalent: \\ \noindent
 -- $y$  is quasi-orthogonal to  $x$; \\ \noindent
 --  for any $\la \in \R$ the vector $x+ \la y$
 lies in the supporting hyperplane to the  ball  $\B_{\norm{x}}(o)$  at $x;$ \\ \noindent
 -- for any $\la \in \R$ the inequality $\norm{x + \la y} \bol \norm{x}$ holds;\\ \noindent
 --  $x$ is orthogonal to  $y$ in the sense of Birkhoff--James (see \cite{DiestelEng}, Ch. 2, \S 1, and \cite{AlonsoMartiniWu_birkhoff_orthogonality}).

Let
\begin{equation*}
\delta_X(\e): = \inf \left\{ 1 - \frac{\|x + y\|}{2}:\ x,y\in\BB,\ \|x -y\| \geqslant \e\right\}
\end{equation*}
and
\begin{equation*}
\mglx{\tau}: = \sup \left\{\frac{\|x + y\| + \|x - y\|}{2} - 1 : \, \|x\| = 1, \|y\| = \tau \right\}.
\end{equation*}
The functions $\delta_X(\cdot): [0,2] \to [0, 1]$ and $\mglx{\cdot}: \R^+ \to \R^+$
 are  referred to as the \emph{moduli of convexity and smoothness of} $X$, respectively.

In \cite{Banass1}  J. Bana{\'s} defined and studied some new modulus of smoothness.
Namely,  he defined
$$\label{mod banass}
\mgbx{\e} = \sup \left\{ 1 - \frac{\norm{x + y}}{2}:\ x,y \in \BB,\ \norm{x -y} \men \e \right\}, \quad \e\in [0,2].
$$

Let $f$ and $g$  be two non-negative functions, each of them defined on a segment $[0, \e].$
We shall say that  $f$ and $g$ are {\it equivalent at zero,} denoted by $f(t) \asymp g(t)$ as $t \to 0,$
if there exist positive constants
$a,b,c,d,e$ such that $a f(bt) \men g(t) \men c f(dt)$ for $t \in [0,e].$

\section{Right-angled triangles}

We will say that a triangle is \emph{right-angled} if one of its legs is quasi-orthogonal to the other one.
(Note that there are completely different ways to define right-angled triangles in normed planes, see also \cite{MartiniAlonsoSpirova_enclosing_disk}.)
In  a Hilbert space this notion coincides with the common, well-known definition of a right-angled triangle.
\begin{remark}
	In  a non-smooth convex Banach space one leg of a triangle can be quasi-orthogonal to the two others.
\end{remark}
For a given right-angled triangle $abc,$ where $ac \prp bc,$ we will say that the
legs $ac, bc$  are the {\it catheti}, and $ab$ the {\it hypotenuse}, of this triangle.

\begin{figure}[ht]%
\center{
\psfrag{x}[1]{\raisebox{-0.9ex}{$\hspace{0.3em}x$}}
\psfrag{o}[1]{\raisebox{0.7ex}{$o$}}
\psfrag{y}[1]{\raisebox{-2ex}{\hspace{0.6em}$y$}}
\psfrag{y1}[1]{\raisebox{-2ex}{\hspace{0.6em}$y_1$}}
\psfrag{y2}[1]{\raisebox{-2.8ex}{$y_2$}}
\psfrag{z}[1]{\raisebox{-1.7ex}{$z$}}
\psfrag{d}[1]{\raisebox{-2.7ex}{$\hspace{-0.4em}d$}}
\psfrag{l}[1]{\raisebox{0.4ex}{\hspace{-0.4em}$l$}}
\includegraphics[scale=0.6]{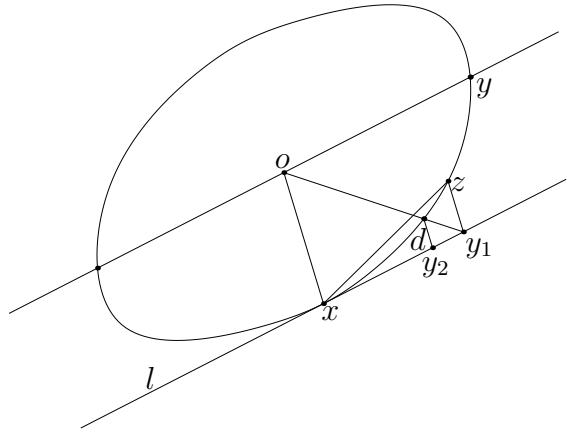}
}
\caption{Right-angled triangles and the unit sphere}
\label{fig_zeta_estimates}
\end{figure}

For  convenience   we draw a simple figure (see Fig. \ref{fig_zeta_estimates}) and  introduce related new moduli by explicit geometric construction.
	Let $x,y \in \SS$ be such that $y \prp x.$ Let $\e \in (0,1],$
	$y_1 = x + \e y.$ Denote by $z$ a point from the unit sphere such that
	$zy_1 \parallel ox$ and  $zy_1 \cap \BB = \{z\}.$  Let $\{d\} = oy_1 \cap \SS.$
	Write $y_2$ for the projection of the point $d$ onto the line $\{x + \tau y: \tau \in \R\}$
	(in the non-strictly convex case we choose $y_2$ such that $dy_2 \parallel ox$).
	Let $p \in J_1(x)$ be such that $\bra p, y \ket = 0$,
i.e., the line $\{x + \tau y: \tau \in \R\}$ lies  in the supporting hyperplane $l = \{a \in X: \bra p, a \ket = 1\}$
of the unit ball at the point $x.$ Then $\norm{zy_1} = \bra p, x-z \ket.$

Consider the right-angled triangle $oxy_1$ (Fig. 1).  In a Hilbert space we have $\norm{oy_1} = \sqrt{1+\e^2},$
but in an arbitrary Banach space the length of the hypotenuse $oy_1$ can vary.
So we introduce moduli that describe the minimal and the maximal length of the hypotenuse in a right-angled triangle in  a Banach space. More precisely, we write
$$
	\zetam{\e}: =  \inf\left\{\norm{x + \e y}: x,y \in \SS, y \prp x\right\}
$$
and
$$
	\zetap{\e}: =  \sup\left\{\norm{x + \e y}: x,y \in \SS, y \prp x\right\}\,,
$$
where $\e$ is an arbitrary positive real number.

In other words, $\zetam{\cdot} - 1$ and $\zetap{\cdot} -1$ describe extrema of the \emph{deviation of a point in a supporting hyperplane from the unit ball}.

On the other hand, the length of the segment $zy_1$ is  the \emph{deviation of a point at the unit sphere from the corresponding
supporting hyperplane}, and at the same time it is a cathetus in the triangle $xzy_1.$

Let $x,y \in \SS$ be such that  $y\urcorner x.$
By definition, put
$$
\la_X(x,y,\e): = \min{\{\la \in \R: \, \norm{x+\e y - \la x} = 1\}}
$$
for any  $\e \in [0,1].$
In the notation of Fig. \ref{fig_zeta_estimates} we have $\la_X(x,y,\e) = \norm{zy_1}.$
The  minimal and the maximal value of $\la_X(x,y,\e)$
characterize the deviation of the unit sphere from an arbitrary supporting hyperplane.
Let us introduce now further moduli.

Define the {\it modulus of supporting convexity} by
$$
\lambda^{-}_{X}(\e) = \inf \{\la_X(x,y,\e):\, x,y \in \BB,\,
 y \urcorner x \}\,,
$$

and the {\it modulus of supporting smoothness} by
$$
\lambda^{+}_{X}(\e) = \sup \{\la_X(x,y,\e):\, x,y \in \BB,\, y \urcorner x \}.
$$

The notions of moduli of supporting convexity and supporting smoothness were introduced and studied  in \cite{Ivanov_supp_modulus}.  These moduli  are very convenient for solving problems concerning the local behaviour of the unit ball compared with that of corresponding supporting hyperplanes.
We will use some of their properties in this paper.

In \cite{Ivanov_supp_modulus} the following inequalities were proved:
\begin{equation}\label{lapx_mglx_ineq}
\mglx{\frac{\e}{2}} \men \lapx{\e} \men \mglx{2\e}, \quad  \e \in\left[0, \frac{1}{2}\right]\,,
\end{equation}
\begin{equation}\label{lamx_mcox_ineq}
     \mcox{\e} \men \lamx{\e}  \men \mcox{2\e}, \quad  \e \in\left[0, 1\right]\,,
\end{equation}
and
\begin{equation}\label{suppconvprop2}
0 \men \lamx{\e} \men  \lapx{\e} \men \e.
\end{equation}

In addition, also a Day-Nordlander type result, referring to these moduli, was proved in \cite{Ivanov_supp_modulus}:
\begin{equation*}
\lamx{\e} \men \la_{H}^{-}(\e) = 1 - \sqrt{1 - \e^2}= \la_{H}^{+}(\e) \men \lapx{\e} \qquad \frl \e \in [0,1].
\end{equation*}

In some sense, moduli of supporting convexity and supporting smoothness are estimates of
a possible value referring to tangents in a Banach space (we fix the length of one of the  catheti and calculate then the minimal and maximal length of the correspondingly other cathetus, which is quasi-orthogonal to the first one).

\begin{remark}\label{remark_la}
	By convexity of the unit ball we have that, for arbitrary $x, y \in \BB$ such that $y \prp x$,
	the function $\la_X(x, y, \cdot)$  is a convex function on the interval $[0, 1].$
\end{remark}

But what can one say about the length of the segment $zy_1$ with fixed norm $\norm{zx}$
(in the notation of the Fig. \ref{fig_zeta_estimates})?

Let us introduce the following new moduli of a Banach space:
\begin{equation}
\phimx{\e} =  \inf \left\{ \bra p, x - z \ket :  x,z\in \SS,\ \|x -z\| \geqslant \e, p \in J_1(x) \right\}\,
\end{equation}
and
\begin{equation}
\phipx{\e} =  \sup \left\{ \bra p, x - z \ket :  x,z\in \SS,\ \|x -z\| \leqslant \e, p \in J_1(x) \right\}\,
\end{equation}
for $\e \in [0,2].$

\begin{remark}
	Due to the convexity of the unit ball we can  substitute inequalities in the definitions of
	$\phimx{\cdot}$ and $\phipx{\cdot}$ to equalities
	(i.e., $\|x -y\| \geqslant \e$  and $\|x -y\| \leqslant \e$   to be $\norm{x-y} = \e$).
\end{remark}

\section{Properties of the catheti}\label{section_catheti}
\begin{lemma}\label{lemma_o_proek_hordi}
In the notation of Fig. \ref{fig_zeta_estimates},
we have $2\norm{y_1x} \bol \norm{xz}.$
\end{lemma}
\begin{prf}
	By the triangle inequality, it suffices to show that $\norm{y_1x} \bol \norm{zy_1}.$
	Let the line $\ell_y$ be parallel to $ox$ with $y \in \ell_y.$
	By construction, we have that the points $x,y_1,z,o,y$ and the line	$\ell_y$
	lie in the same plane  -- the linear span of the vectors $x$ and $y$.
	So the lines $\ell_y$ and $xy_1$ intersect, and by $c$ we denote their intersection point.
	Note that  $oycx$ is a parallelogram and $\norm{yc} = 1$;
	the segment $yx$ belongs to the unit ball and does not intersect the interior of the segment $zy_1.$
	Let $\{z'\} = zy_1 \cap yx.$
 	By similarity, we have
	$$
		\norm{zy_1} \men \norm{y_1 z'} = \frac{\norm{xy_1}}{\norm{x  c}} \norm{y  c} = \norm{x y_1}.
	$$
\end{prf}
\begin{figure}[h]%
\center{
\psfrag{x}[1]{\raisebox{0.2ex}{$\hspace{0.3em}x$}}
\psfrag{o}[1]{\raisebox{0.3ex}{$o$}}
\psfrag{c}[1]{\hspace{0.3em}$c$}
\psfrag{y}[1]{\raisebox{-1.6ex}{\hspace{0.5em}$y$}}
\psfrag{y1}[1]{\raisebox{-2ex}{\hspace{0.3em}$y_1$}}
\psfrag{zz}[1]{\raisebox{1ex}{\hspace{-0.8em}$z'$}}
\psfrag{z}[1]{\raisebox{-0.4ex}{\hspace{0.2em}$z$}}
\psfrag{l}[1]{\hspace{0.6em}$\ell_y$}
\psfrag{ll}[1]{\raisebox{0.9ex}{\hspace{0.0em}$l$}}
\includegraphics[scale=0.5]{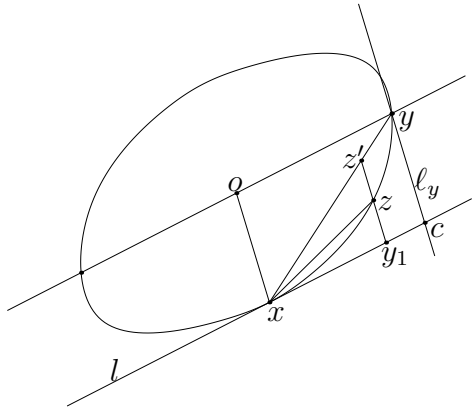}
}
\caption{ Under the conditions of Lemma \ref{lemma_o_proek_hordi} we have $2\norm{xy_1} \bol \norm{xz}.$}
\end{figure}
It is worth noticing that under the conditions of   Lemma \ref{lemma_o_proek_hordi} we have
	that $y_1$ is a projection along the vector $ox$ of the point $z$
	 on some supporting hyperplane of the unit ball at $x$.
	Moreover, $y_1$  belongs to the metric projection of the point $y$ on this hyperplane.
	In other words, Lemma \ref{lemma_o_proek_hordi} shows us that
	if one projects the segment $xz$ along the vector $ox$ onto the hyperplane
	which supports the unit ball at $x,$
	then the length of the segment decreases no more than by a factor of 2.

\begin{lemma}\label{lemma_phimpx_lampx}
 Let  $X$ be an arbitrary Banach space. Then the following inequalities  hold:
 \begin{equation}\label{phimx_inequality}
 \lamx{\frac{\e}{2}} \leqslant \phimx{\e} \leqslant \lamx{2\e}\, \mbox{ and }
\end{equation}
 \begin{equation}\label{phipx_inequality}
  \lapx{\frac{\e}{2}} \leqslant \phipx{\e} \leqslant \lapx{2\e}
 \end{equation}
 for $\e \in [0,1/2].$
\end{lemma}
\begin{prf}
In the notation of Fig.\ref{fig_zeta_estimates} we assume that for arbitrary $x, y$ with $y \prp x$  the equality
$\norm{zx} = \e$ holds.  Then $\la_X(x, y, \norm{xy_1}) = \norm{y_1 z}.$
Let $p \in J_1(x)$ be such that $\langle p, y\rangle = 0.$ Hence $\norm{y_1z} = \bra p, x- y\ket.$
Since $\norm{xy_1} \men \norm{y_1z} + \norm{zx} \men 2\e,$ and
taking into account Lemma \ref{lemma_o_proek_hordi},   we get
$$
	 \frac{\e}{2} \leqslant \norm{xy_1} \leqslant 2\e \leqslant 1.
$$
Due to this and by Remark \ref{remark_la}, we have
$$
	\la_X\!\! \left(x, y, \frac{\e}{2}\right) \leqslant \bra p, x - y  \ket \leqslant  \la_X\left(x, y, 2\e \right).
$$
	Taking infimum (supremum) on the right-hand side, left-hand side or in the middle part of
	the last inequality, we obtain \reff{phimx_inequality} and \reff{phipx_inequality}.
\end{prf}

From Lemma \ref{lemma_phimpx_lampx} and the inequalities \reff{lamx_mcox_ineq} and \reff{lapx_mglx_ineq} we have the following corollary.
\begin{corollary}
	Let  $X$ be an arbitrary Banach space. Then  $\phipx{\e}\asymp\mglx{\e}$  and
	 $\phimx{\e} \asymp \mcox{\e}$ as $\e \to 0$, and for $\e \in\left[0, \frac{1}{2}\right]$
	the following inequalities  hold:
$$
 \mglx{\frac{\e}{4}} \men \phipx{\e} \men \mglx{4\e},  \, \mbox{ and }
$$
$$
     \mcox{\e} \men \phimx{\e}  \men \mcox{4\e} .
$$
\end{corollary}

Now we will prove a Day-Nordlander type result for $\phimx{\cdot}$ and $\phipx{\cdot}.$
Let us suitably generalize the notion of modulus of convexity and the notion of Bana{\'s} modulus. Namely, let
\begin{equation*}
\mcoxt{\e}{t} = \inf \left\{ 1 - \frac{\|t x + (1-t)y\|}{2}:\ x,y\in\SS,\ \|x -y\| = \e\right\}
\end{equation*}
and
\begin{equation*}
\mgbxt{\e}{t} = \sup \left\{ 1 - \frac{\|t x + (1-t)y\|}{2}:\ x,y\in\SS,\ \|x -y\| = \e\right\}\,,
\end{equation*}

respectively. Using the same method as in the classical paper \cite{Day-Nord}, we get
\begin{lemma}
	Let  $X$ be an arbitrary Banach space. Then the following inequalities hold:
	\begin{equation}
	\mcoxt{\e}{t} \leqslant \delta_H (\e, t) = 1 - \sqrt{1 - t(1-t)\e^2} = \delta_{H}^{+} (\e, t) \leqslant \mgbxt{\e}{t}.
	\end{equation}
\end{lemma}
\begin{prf}
Since the proof is almost the same as in \cite{Day-Nord}, we present only a short sketch. Clearly, again
it is sufficient to prove the lemma in the two-dimensional case.

If the two unit vectors $x = (x_1, x_2)$ and $y = (y_1, y_2)$ are rotated around the unit circle,
while their difference $x-y$ has constantly the norm $\e,$
the endpoint of the vector $t x + (1-t) y$ describes a curve $\Gamma_t.$

The following integral expresses the area of the region inside the curve described by the endpoint of the vector $x - y,$ if this vector is laid off from a fixed point:
$$
	\int (y_1 - x_1) d(y_2 - x_2)\,.
$$
On the other hand, the mentioned curve is a homothet of the unit circle with ratio
$\e.$ Hence this integral equals $\e^2 A,$ where $A$ is the area of the unit ball ($A = \int x_1 d x_2 = \int y_1 d y_2$).
From this we have that
$$
\int x_1 d y_2 + \int y_1 d x_2 = 2A - \e^2A.
$$
Now it is clear that the area of the region inside $\Gamma_t$ equals
$$
	\int (t x_1 + (1-t) y_1) d (t x_2 + (1-t) y_2) = A (1 - t(1-t)\e^2).
$$
Hence continuity arguments imply that there exists a point $z \in \Gamma_t$ with the norm $\sqrt{1 - t(1-t)\e^2}.$
\end{prf}
\begin{theorem}\label{theorem_daynord_phi}
   	Let  $X$ be an arbitrary Banach space. Then the following inequalities  hold:
	\begin{equation}
	\phimx{\e} \leqslant \varphi_H^{-}(\e) = \frac{\e^2}{2} =  \varphi_H^{+}(\e) \leqslant \phipx{\e}.
	\end{equation}
\end{theorem}
\begin{prf}
It is sufficient to prove the theorem in the two-dimensional case.
Let $x \in \SS$ and $p \in J_1(x).$

Assume that $X$ is a uniformly smooth space.
Notice that $p$ is a Frechet derivative of the norm at the point $x.$
Taking into account that $\BB$ is convex, for an arbitrary $y$ we have
$$
	\bra p, x - y \ket =
	\lim\limits_{t \searrow 0} \frac{\norm{x} -  \norm{x + t(y-x)}}{t} = \liminf\limits_{t > 0}
	\frac{1 -  \norm{x + t(y-x)}}{t}.
$$
Fix an arbitrary $\gamma > 0.$
Since $X$ is uniformly smooth, there exists a $t_0 < \gamma$ such that
for  arbitrary $x, y \in \SS$, $ \norm{y-x} = \e$ and $t \in (0, t_0)$
we have
\begin{equation}\label{Th_Day_Nord_phi_1}
\frac{1 -  \norm{x + t(y-x)}}{t} - \gamma	\leqslant \bra p, x - y \ket \leqslant \frac{1 -  \norm{x + t(y-x)}}{t}.
\end{equation}

Taking the infimum (supremum) in the last line, we get
$$
\frac{\mcoxt{\e}{t}}{t}	 - \gamma \leqslant \phimx{\e} \leqslant \frac{\mcoxt{\e}{t}}{t}
$$
$$
\left(
	\frac{\mgbxt{\e}{t}}{t}- \gamma \leqslant \phipx{\e} \leqslant \frac{\mgbxt{\e}{t}}{t}
\right).
$$
Passing to the limit as $\gamma \to 0,$
we have
$$
\phimx{\e}  =  \lim_{t \to 0}\frac{\mcoxt{\e}{t}}{t} \leqslant
\lim_{t \to 0}\frac{\delta_H (\e, t)}{t} = \frac{\e^2}{2}
$$
$$
\left(
\phipx{\e}  =
 \lim_{t \to 0}\frac{\mgbxt{\e}{t}}{t} \geqslant
\lim_{t \to 0}\frac{\delta_H^{+}(\e, t)}{t} = \frac{\e^2}{2}
\right).
$$

Let us now consider the case of a non-smooth space $X.$
  Let $SP$ be the set of all points of smoothness at the unit circle.
We know that the unit circle is compact.
Then there exists $t_0 < \gamma$ such that
for  arbitrary  $x \in SP$, $y \in \SS$, $ \norm{y-x} = \e$ and $t \in (0, t_0)$
we can write the inequality \reff{Th_Day_Nord_phi_1}.

Moreover, the set $\SS \setminus SP$ has  measure zero. Thus, the infimum (supremum)
of
$
{1 -  \norm{x + t(y-x)}}
$
taken over all $x \in SP$ coincides with $\mcoxt{\e}{t}$ ($\mgbxt{\e}{t}$).
So we have
$$
\phimx{\e}  \leqslant  \limsup_{t \to 0}\frac{\mcoxt{\e}{t}}{t} \leqslant \frac{\e^2}{2}
$$
and
$$
\phipx{\e}  \geqslant
 \liminf_{t \to 0}\frac{\mgbxt{\e}{t}}{t} \geqslant \frac{\e^2}{2}.
$$
\end{prf}

\section{Properties of the Hypotenuse}\label{section_hypothenuse}
\begin{lemma}
	Let  $X$ be an arbitrary Banach space. Then for $ \e \in [0, 1]$ the following inequalities  hold:
	\begin{equation} \label{zetam_inequality}
		\lamx{\frac{\e}{1+\e}} \leqslant \zetam{\e} -1 \leqslant \lamx{\e}\,,
	\end{equation}
	\begin{equation} \label{zetap_inequality}
		\lapx{\frac{\e}{1+\e}} \leqslant \zetap{\e} - 1 \leqslant \lapx{\e}\,.
	\end{equation}
\end{lemma}
\begin{prf}
	From the triangle inequality we have that $\norm{y_1d}$ equals the distance from the point $y_1$
	to the unit ball. Hence
	\begin{equation}\label{zeta_estimates_1}
	\norm{y_1d} \leqslant \norm{y_1z} = \lambda_X\left(x,y,\e \right) \leqslant \e.
	\end{equation}
	By similarity arguments and \reff{zeta_estimates_1} we have
	$$
	\norm{xy_2} = \frac{\norm{od}}{\norm{od} + \norm{dy_1}} \norm{xy_1}
	= \frac{1}{1 +\norm{dy_1}}\e \geqslant
	\frac{\e}{1 + \e}.
	$$
	Then, by construction and by the convexity of the unit ball, we get
	the inequality	
	\begin{equation}\label{zeta_estimates_2}
			\norm{y_2d} = \lambda_X(x, y, \norm{xy_2}) \geqslant \lambda_X\!\!\left(x, y, \frac{\e}{1 + \e} \right).
	\end{equation}
	Since $y_2$ is a projection of the point $d$ onto the line $\{x + \tau y: \tau \in \R\},$ we have
	$\norm{y_2d} \leqslant \norm{dy_1}.$ Combining the previous inequality with \reff{zeta_estimates_1} and
	(\ref{zeta_estimates_2}), we obtain the inequalities
	$$	
		\lambda_X\!\!\left(x, y, \frac{\e}{1 + \e} \right)  \leqslant \norm{dy_1} \leqslant
		 \lambda_X\!\left(x,y,\e \right).
	$$
	Taking infimum (supremum) on the right-hand side, left-hand side or in the middle part of
	the last line, we obtain \reff{zetam_inequality} and \reff{zetap_inequality}.
\end{prf}
\begin{corollary}
	Let  $X$ be an arbitrary Banach space. Then $\zetap{\e} - 1\asymp\mglx{\e}$  and
	 $\zetam{\e} - 1 \asymp \mcox{\e}$ as $\e \to 0$, and
	the following inequalities  hold:
$$
 \mglx{\frac{\e}{2(1+\e)}} \men \zetap{\e} \men \mglx{2\e}, \quad  \e \in\left[0, \frac{1}{2}\right]\,, \mbox{ and }
$$
$$
     \mcox{\frac{\e}{1+\e}} \men \zetam{\e}  \men \mcox{2\e}, \quad  \e \in\left[0, 1\right].
$$
\end{corollary}

Now we will prove results of Day-Nordlander type for $\zetam{\cdot}$ and $\zetap{\cdot}.$

Suppose we have an orientation $\omega$ in $\R^2.$
We will say that a curve $C$ in the plane is a {\it good curve} if
it is a closed rectifiable simple Jordan curve, which is enclosed by a star-shaped set $S$ with  center at the origin
and continuous radial function.
\begin{lemma}\label{lemma_integration}
Let $C_1$  be a closed  simple Jordan curve enclosing
the convex set $S_1$ with area $A_1 > 0$ and $0 \in \vn S_1.$ Let $C_2$ be a good curve, which is enclosing an area of measure $A_2$. Then
\begin{enumerate}
\item we can parametrize $C_i$  by a function $f^i(\cdot):[0,1) \to C_i$ ($i =1,2$) in such a way that
\begin{enumerate}
\item $f^2(\tau)$  is a direction vector of the supporting line
of the set $S_1$ at the point $f^1(\tau)$ for all $\tau \in [0,1);$
\item $[f^1(\tau)  , f^2(\tau)] = \omega$ for all $\tau \in [0,1);$
\item the functions $f^i(\cdot)$ ($i=1,2$) are angle-monotone;
\end{enumerate}
\item the curve $C_3 = \left\{f^1(\tau)+f^2(\tau): \tau \in [0, 1)\right\}$ encloses an area of measure $A_1 + A_2.$  	
\end{enumerate}
\end{lemma}
\begin{prf}
1) First of all, due to the continuity of the radial function of the curve $C_2$
we can assume that $C_2$ and $C_1$ are coincident.

Let  $C_1$ be a smooth curve. Let $f^1: [0,1) \to C_1$ be a parametrization given by clockwise rotation.
Then at every point $f^1(\tau)$ we have a unique supporting line to $S_1$, and we can choose $f^2(\tau)$
in a proper way.   In this case the problem is quite  easy and one can see its geometric interpretation.

The general case (when $C_1$ has non-smooth points) yields additional difficulties.
At a point of non-smoothness we have continuously many supporting lines; hence we cannot give a parametrization  depending only on this point of $C_1.$
However, in \cite{joly_parametriztion} Joly gives a suitable parametrization.\\
\noindent
2) Let $A_3$ be the measure of the area enclosed by $C_3.$
Let $f^i(\cdot)$  be the parametrization of $C_i$ ($i=1,2$) constructed above.  Fix $\mu \in \R.$
Denote by $S(\mu)$ and $A(\mu)$ the set and the area enclosed
by the curve $C(\mu) = \left\{f^1(\tau)+ \mu f^2(\tau): \tau \in [0, 1)\right\}$, respectively.
Since for all $\tau \in [0,1)$ we have that  $f^2(\tau)$  is a direction vector of the supporting line
of the set $S_1$ at the point $f^1(\tau),$ then we have $S_1 \subset S(\mu).$
Hence $A(\mu) \geqslant  A_1.$ Using consequences of Green's formula and
properties of the Stieltjes integral, we have
$$
\int\limits_{\tau \in [0,1)} f^1_1 d f^1_2 \leqslant 	
\int\limits_{\tau \in [0,1)} (f^1_1(\tau)+ \mu f^2_1(\tau)) d(f^1_2(\tau)+ \mu f^2_2(\tau)).
$$
Therefore, for all $\mu \in \R$ the following inequality holds:
\begin{equation*}
	\mu^2 \int\limits_{\tau \in [0,1)} f^2_1 d f^2_2 +
	\mu \left(\int\limits_{\tau \in [0,1)} f^1_1 d f^2_2 + \int\limits_{\tau \in [0,1)} f^2_1 d f^1_2 \right) \geqslant 0.
\end{equation*}
This implies that
$$
\left(\int\limits_{\tau \in [0,1)} f^1_1 d f^2_2 + \int\limits_{\tau \in [0,1)} f^2_1 d f^1_2 \right) = 0.
$$
So we have
$$
	A_3 =  A(1) = \int\limits_{\tau \in [0,1)} (f^1_1(\tau)+  f^2_1(\tau)) d(f^1_2(\tau)+  f^2_2(\tau)) =
	\int\limits_{\tau \in [0,1)} f^1_1 d f^1_2  + \int\limits_{\tau \in [0,1)} f^2_1 d f^2_2 = A_1+ A_2.
$$
\end{prf}

\begin{theorem}
   	Let  $X$ be an arbitrary Banach space. Then the following inequalities hold:
	\begin{equation}\label{daynordzeta}
	\zetam{\e} \leqslant \zeta_H^{-}(\e) = \sqrt{1 + \e^2} =  \zeta_H^{+}(\e) \leqslant \zetap{\e}.
	\end{equation}
\end{theorem}
\begin{prf}
 Again it is sufficient to prove the theorem in the two-dimensional case.
 Applying Lemma \ref{lemma_integration} for $C_1 = \SS, C_2 = \partial \B_\e (o)$ and using continuity arguments
 we obtain \reff{daynordzeta}.
\end{prf}
\begin{remark}
In \cite{joly_parametriztion} inequality \reff{daynordzeta} was proved for the subcase $\e = 1.$
\end{remark}
\section{Some notes about monotonicity properties of the dual mapping} \label{section_monotonicity}
The notion of monotone operator is well-known and has a lot of applications and useful generalizations. 
Let us recall some related notions and, based on them, explain their relations to the geometry of the unit sphere.

Let $X$ be a Banach space, $T: X \to X^*$ a point-to-set operator, and $G(T)$ its graph.
Suppose that  the following inequality holds:
\begin{equation}\label{inequality_monotonicity}
\bra  p_x - p_y, x - y \ket \geqslant \alpha \norm{x - y}^2  \quad \mbox{ for all } (x, p_x), (y, p_y) \in G(T).
\end{equation}
If 
\begin{enumerate}
\item
$\alpha = 0,$ then $T$ is a monotone operator. For example, the subdifferential of a convex function is a monotone operator. 
\item $\alpha > 0 ,$ then $T$ is a strongly monotone operator. For example, the subdifferential of a strongly convex function on a Hilbert space is a strongly monotone operator. 
\item $\alpha < 0 ,$ then $T$ is a hypomonotone operator. For example, 
the subdifferential of a prox-regular function on a Hilbert space  
is a hypomonotone operator (see \cite{IvMonEng}). 
\end{enumerate} 

Inequality \reff{inequality_monotonicity} is often called the {\it variational inequality.}
Usually, the operator $T$ is a derivative or subderivative of a convex function.
So we can speak about the variational inequality for a convex function.

As usual in  convex analysis,  we can reformulate inequality \reff{inequality_monotonicity} for convex
(or prox-regular) sets and their normal cone (or Frechet normal cone) (see \cite{Pol_Rock_Thi1}), in this case $T(x)$ is a intersection of the $\SSS$ and the normal
cone to the set at point $x$.  
In a Hilbert space there are some characterizations of strongly convex and prox-regular functions  
(or strongly convex and prox-regular sets) via  the variational inequality (see \cite{Clark}, \cite{proxreg_rockafellar} and \cite{Pol_Rock_Thi1}). 

But in a Banach space the situation is much more complicated and it is getting obvious that the right-hand  side of the variational inequality cannot always be a quadratic function. So, in many applications we have to 
substitute $\alpha \norm{x - y}^2$ in \reff{inequality_monotonicity} by some proper convex function 
$\alpha(\norm{x- y}).$

For example, what can we say about the  most simple convex function in a Banach space -- its norm 
(in this case $T$ is a dual mapping)?
Even in a Hilbert space, for  arbitrary $x, y$ we can only put  zero in the right-side 
of the variational inequality. Nevertheless, there exist variational inequalities for norms depending on 
$\norm{x}, \norm{y},$ and $\norm{x-y}$. For example, in \cite{Charinequalities_uc_usm} characterizations 
of  uniformly smooth and uniformly convex Banach spaces were given in terms of monotonicity properties of the dual mapping.

In this paragraph we investigate monotonicity properties of the dual mapping onto the unit sphere.
In fact, we study   monotonicity properties of the convex function  on its Lebesgue level.
Hence this results can be generalized to an arbitrary convex function.  

We are interested in  asymptotically tight lower and upper bounds for the 
$\bra p_1 - p_2, x_1 - x_2 \ket,$ where  $x_1, x_2 \in \SS,  p_1 \in J_1 (x_1), p_2 \in J_1(x_2).$ 
For the sake of convenience we introduce new moduli:
$$
	\gammap{\e} = \sup\{\bra p_1 - p_2, x_1 - x_2 \ket : x_1, x_2 \in \SS, \norm{x_1 - x_2} = \e, p_1 \in J_1 (x_1), p_2 \in J_1(x_2) \}
$$
and
$$
	\gammam{\e} = \inf\{\bra p_1 - p_2, x_1 - x_2 \ket : x_1, x_2 \in \SS, \norm{x_1 - x_2} = \e, p_1 \in J_1 (x_1), p_2 \in J_1(x_2) \}
$$
for each $\e \in [0, 2].$

\begin{figure}[ht]%
\center{
\psfrag{x}[1]{\raisebox{0.3ex}{$\hspace{0.3em}x$}}
\psfrag{x1}[1]{\raisebox{-1.4ex}{$\hspace{0.3em}x_1$}}
\psfrag{x2}[1]{\raisebox{0.7ex}{$\hspace{0.9em}x_2$}}
\psfrag{o}[1]{\raisebox{0.1ex}{$o$}}
\psfrag{y}[1]{\raisebox{0.1ex}{\hspace{0.15em}$y$}}
\psfrag{z1}[1]{\raisebox{1.0ex}{\hspace{0.6em}$z_1$}}
\psfrag{z2}[1]{\raisebox{1.4ex}{$\hspace{-1.4em}z_2$}}
\psfrag{l}[1]{\raisebox{-2.8ex}{$\hspace{-1.2em}l$}}
\psfrag{l1}[1]{\raisebox{-1.9ex}{\hspace{-0.4em}$l_1$}}
\psfrag{l2}[1]{\raisebox{0.4ex}{\hspace{0.em}$l_2$}}
\psfrag{b2}[1]{\raisebox{0.4ex}{\hspace{-1.3em}$b_2$}}
\psfrag{b1}[1]{\raisebox{-2.4ex}{\hspace{-2.em}$b_1$}}
\includegraphics[scale=0.6]{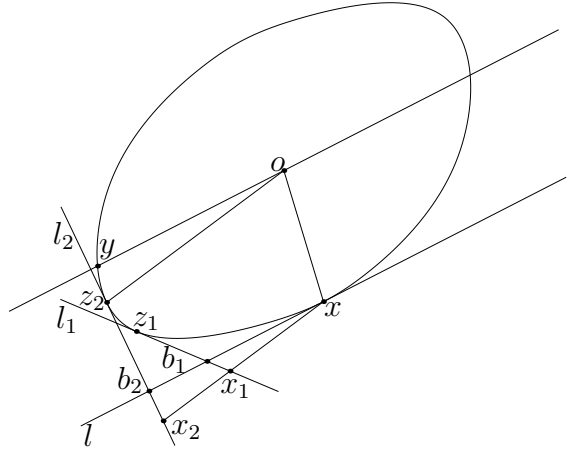}
}
\caption{Illustration of the proof of Lemma  \ref{lemma_gamma_1}}
\label{fig_gamma_monot}
\end{figure}

\begin{lemma}\label{lemma_gamma_1}
Let  $X$ be an arbitrary Banach space. Then the functions $\gammap{\cdot}$ and $\gammam{\cdot}$ are monotonically increasing functions on $[0, 2]$.
\end{lemma}
\begin{prf}
In the notation of Fig. \ref{fig_zeta_estimates}, let $z_1, z_2$ be points in the arc $-xyx$ of the unit circle such that $z_1$ belongs to the arc $xz_2$ (here and in the sequel all arcs lie in the plane of $xoy$).
Let $p \in J_1(x), q_1 \in J_1(z_1), q_2 \in J_1(z_2).$ 
It is worth mentioning that $\norm{xz_1} \leqslant \norm{xz_2}$ (see \cite{AlonsoBenitez}, Lemma 1). 
So, to prove our Lemma it is sufficient to show that 
\begin{equation}\label{gamma_monot_1}
 \bra p  - q_1,  x - z_1 \ket \leqslant \bra  p  - q_2,  x - z_2 \ket
\end{equation} 
From the convexity of the unit ball we have that $\bra  p , x - z_1 \ket \leqslant \bra  p,  x - z_2 \ket.$ 
To prove inequality \reff{gamma_monot_1}, let us show that $\bra q_1,  z_1 - x \ket \leqslant \bra  q_2,  z_2 - x \ket.$ 

We can assume that $X$ is the plane of $xoy.$ 
By definition, put $l = \{a \in X: \bra a, p \ket = 1 \}$,
$l_1 = \{a \in X: \bra a, q_1 \ket = 1 \},$  $l_2 = \{a \in X: \bra a, q_2 \ket = 1 \},$ and
$H^+ = \{p \in X: \bra a, p \ket \geqslant 1 \}.$  

The first case: let $z_2$ be in the arc $xy$ of the unit circle (see Fig. \ref{fig_gamma_monot}).
All three cases $l = l_1,$ $l = l_2$ or $l_1 = l_2$ are trivial.
Let $l \cap l_1 = \{b_1\},$ $l \cap l_2 = \{b_2\}.$ 
Again, all three cases $x = b_1,$ $x = b_2$ or $b_1 = b_2$ are trivial.
By  convexity arguments, $b_1$ 
belongs to the relative interior of the segment $xb_2$ and $l_1 \cap l_2 \notin H^+.$ 
Hence $l_1$ separates point $x$ and the ray $l_2 \cap H^+$ in the  half-plane $H^+.$ 
Let $x_2$ be a projection of the point $x$ onto $l_2$
(in the non-strictly convex case we choose $x_2$ such that $xx_2 \parallel oz_2$).
Then the segment $xx_2$ is parallel to $oz_2,$ and therefore $xx_2 \subset H^+.$ 
Now we can say that the segment $x_2x$ and the line $l_1$ have an intersection point; let it be  $x_1.$
Since the values $\bra q_1,  z_1 - x \ket$ and $\bra  q_2,  z_2 - x \ket$ are equal to the distances from the  point 
$x$ to the lines $l_1$ and $l_2$, respectively, we have:
$$
\bra q_1,  z_1 - x \ket \leqslant \norm{xx_1} <  \norm{xx_2} = \bra  q_2,  z_2 - x \ket.
$$

The second case: let $z_2$ be in the arc $-xy$ of the unit circle.
We can assume that $z_1$ lies on the arc $-xy$ of the unit circle, too 
(if $z_1$ lies on the arc $xy$ of the unit circle, 
by the first case we can substitute $z_1$ to $y$).
We have that $\bra -q_i, -z_i - x \ket = 2 - \bra q_1, z_i - x \ket$ for $i = 1,2.$
Therefore, applying the first case to the points $-z_1, -z_2, x $ and
to the functionals $p, -q_1, -q_2,$  we have proved  the second case.
\end{prf}
\begin{remark}
	It is worth mentioning that in the first case of Lemma \ref{lemma_gamma_1} the lines $l_1$ and $ox$ 
	can have no common point in  $H^+.$
\end{remark}
\begin{remark}
Using Lemma \ref{lemma_gamma_1}, we can modify the definitions of $\gammap{\cdot}$ and $\gammam{\cdot}$ by
$$
	\gammap{\e} = \sup\{\bra p_1 - p_2, x_1 - x_2 \ket : x_1, x_2 \in \SS, \norm{x_1 - x_2}  \leqslant \e, p_1 \in J_1 (x_1), p_2 \in J_1(x_2) \}
$$
and
$$
	\gammam{\e} = \inf\{\bra p_1 - p_2, x_1 - x_2 \ket : x_1, x_2 \in \SS, \norm{x_1 - x_2} \geqslant \e, p_1 \in J_1 (x_1), p_2 \in J_1(x_2) \}
$$
for each $\e \in [0, 2].$
\end{remark}
 \begin{lemma}
	Let  $X$ be an arbitrary Banach space. 
	Then  the following inequalities  hold:
	\begin{equation} \label{gammap_inequality}
		\phipx{\e} \leqslant \gammap{\e} \leqslant 2 \phipx{\e} \, \mbox{ for } \e \in [0,2],
	\end{equation}
	\begin{equation} \label{gammam_inequality}
		2\phimx{\frac{e}{4}} \leqslant \gammam{\frac{\e}{4}} \leqslant  \phimx{\e} \mbox{ for } \e \in [0,1]\,.
	\end{equation}
\end{lemma}
\begin{prf}
	All inequalities, except for the right-hand side of  \reff{gammam_inequality}, are obvious.
	
	Let us prove that $\gammam{\frac{\e}{4}} \leqslant 2 \phimx{\e}.$
	It is sufficient to prove the lemma in the two-dimensional case.
	In this case and in the notation of Fig. \ref{fig_zeta_estimates} we can put 
	$\norm{zx} = \e$ and $\norm{y_1z} = \phimx{\e}$.
	Let $y_b$ be a bisecting point of the segment $xy_1.$ Denote by $z_b$ a point from the unit sphere such that
	$z_by_b \parallel ox$ and  $z_by_b \cap \BB = \{z_b\}.$  
	Let $p_b \in J_1(z_b).$ Denote by $l_b$ the line $\{a \in X : \bra p_b, a\ket  = 1\}.$
	By convexity the line $l_b$ intersects the segment $zy_1,$ and we denote the intersection point as $a_1.$
	By definition put $\{a_2\} = l_1 \cap \{\tau x : \tau \in \R\}.$ 
	From the trapezoid $a_2xa_1y_1$ we have that
	\begin{equation}\label{lemma_gamma_2_1}
	\norm{y_bz_b} + \norm{xa_2} \leqslant \norm{y_1a_1} \leqslant \norm{zy_1} = \phimx{\e}.
	\end{equation}
	Since $\bra p_b, z_b - x\ket $ equals  the distance from the point $x$ to the line $l_b,$
	we have that $\bra p_b, z_b - x\ket \leqslant \norm{xa_2}.$ 	
	From here, since $\bra p, x- z_b\ket = \norm{y_bz_b},$ and  from inequality \reff{lemma_gamma_2_1} we obtain 
	\begin{equation*}
	 \bra p - p_b, x- z_b \ket \leqslant \phimx{\e}.
 	\end{equation*}
 	From Lemma \ref{lemma_gamma_1} it is sufficient to show that $\norm{xz_b} \geqslant \frac{\e}{4}.$
 	By definition put $\{z'\} = y_bz_b \cap xz.$
	Obviously, we have that 
	$$
	\norm{xz_b} \geqslant \norm{xz'} - \norm{z'z_b} \geqslant \norm{xz'} - \norm{z'y_b} = \frac{\e - \phimx{\e}}{2}. 
	$$ 	
 	Using Theorem \ref{theorem_daynord_phi}, we see that 
 	$$
 		\norm{xz_b} \geqslant \frac{\e}{2} - \frac{\e^2}{4} \geqslant \frac{\e}{4}.
 	$$
\end{prf}

\begin{corollary}
	Let  $X$ be an arbitrary Banach space. Then  $\gammap{\e}\asymp\mglx{\e}$  and
	 $\gammam{\e} \asymp \mcox{\e}$ as $\e \to 0$ and for $\e \in\left[0, \frac{1}{2} \right]$
	the following inequalities  hold:
$$
 \mglx{\frac{\e}{4}} \men \gammap{\e} \men 2\mglx{4\e}  \, \mbox{ and }
$$
$$
     2\mcox{\frac{\e}{4}} \men \gammam{\frac{\e}{4}}  \men \mcox{\e} .
$$
\end{corollary}
\begin{remark}
Combining results from  \cite{Charinequalities_uc_usm} for some constant $c_1, c_2, c_3, c_4$ (depending on $X$) one can get the following inequality:
$$
	c_1 \mglx{c_2\e} \geqslant \gammap{\e} \geqslant \gammam{\e} \geqslant c_3 \mcox{c_4 \e}. 
$$
\end{remark}

\section{Some open questions}\label{section_open_questions}
Although there are no difficulties to prove an analogue of  the Day-Nordlander theorem for the moduli $\gammap{\cdot}$ and $\gammam{\cdot}$ moduli  in the infinite-dimensional case
using Dvoretzky's theorem (see \cite{Dvoretzky}),
we have no proof  for the following conjecture in the finite-dimensional case:
\begin{conjecture}	\label{conjecture_DayNord_gamma}
 	Let  $X$ be an arbitrary Banach space. Then the following inequalities hold:
	\begin{equation}
	\gammam{\e} \leqslant \gamma^-_H(\e) = \e^2 = \gamma^+_H(\e) \leqslant \gammap{\e}.
	\end{equation}
\end{conjecture}

All moduli mentioned above characterize certain geometrical properties of the unit ball.
Obviously, the geometry of the unit ball totally describes the geometry of the unit ball in the dual space.
Nevertheless, we know a few results about coincidences of values of some moduli or other characteristics of a Banach space and its dual space. We are interested in properties of the dual mapping (i.e., $x \to J_1(x)$).
The following conjecture seems to be very essential.
By definition, put
$$
	d^-_X(\e) =
	\inf\{\norm{p_1 - p_2}| p_1 \in J_1(x_1), p_2 \in J_1(x_2), \norm{x_1 - x_2} = \e, x_1, x_2 \in \SS \}
$$
and
$$
	d^+_X(\e) =
	\sup\{\norm{p_1 - p_2}| p_1 \in J_1(x_1), p_2 \in J_1(x_2), \norm{x_1 - x_2} = \e, x_1, x_2 \in \SS \}.
$$
\begin{conjecture}
 	Let  $X$ be an arbitrary Banach space. Then the following inequalities hold:
	\begin{equation}
	d^-_X(\e) \leqslant d^-_H(\e) = \e = d^+_H(\e) \leqslant d^+_X(\e).
	\end{equation}
\end{conjecture}

It is well-known that the equality $\mcox{\e} = \delta_H(\e)$ for $\e \in [0,2)$ implies that 
$X$ is an inner product space (see \cite{Day_innerproduct}). There exist such results for some other moduli
(See \cite{AlonsoBenitez} and \cite{amir_char_inner_pr_sp}). We are interested in the following question:
\begin{question}
For what modulus $f_X(\cdot)$ 
$(f_X(\cdot) = \phimx{\cdot}, \phipx{\cdot}, \zetam{\cdot}, \zetap{\cdot}, \lamx{\cdot}, \lapx{\cdot})$ does the equality $f_X(\e) = f_H(\e),$ holding for all $\e$ in the domain of the function $f_X(\cdot)$ (or even for fixed $\e$), imply that $X$ is an
inner product space?
\end{question}

The definitions of the  moduli $\zetam{\cdot} - 1$ and   $\zetap{\cdot} - 1$ are similar to the definitions of  Milman's moduli,
which were introduced in \cite{Milman_eng}  as
$$
 \beta^-_X(\e)=  \inf\limits_{x,y \in \SS}\{\max\{\norm{x + \e y}, \norm{x - \e y}\}-1\}
$$
and
$$
 \beta^+_X(\e)=  \sup\limits_{x,y \in \SS}\{\min\{\norm{x + \e y}, \norm{x - \e y}\}-1\}.
$$

We think that in the definitions of  Milman's moduli it is sufficient to take only $y \prp x.$ 
Hence we get
\begin{conjecture}
Let  $X$ be an arbitrary Banach space. Then for positive $\e$ we have 
$$
	\zetam{\e} - 1 =\beta^-_X(\e) \mbox{  and  } \zetap{\e} - 1 = \beta^+_X(\e).
$$
\end{conjecture}

\nocite{Banas2}
\nocite{modulsmooth1}
\nocite{Guirao_2}
\nocite{banas1990convexity}
\nocite{AlonsoBenitez}
\nocite{aboutMilmansmoduli}
\nocite{Serb_DayNord}

\end{document}